\documentclass[11pt]{article} 
\usepackage{hyperref,amsmath,amsthm,amsfonts,amssymb,amscd,enumerate}
 
\theoremstyle{plain} 
\newtheorem{theorem}{Theorem}

\newtheorem{corollary}[theorem]{Corollary} 
\newtheorem{lemma}[theorem]{Lemma} %

\theoremstyle{defintion}
\newtheorem{remark}[theorem]{Remark}

\newcommand{\ca}{\mathcal {A}}

\newcommand{\ch}{\mathcal {H}}

\newcommand{\cu}{\mathcal {U}}

\newcommand{\bc}{{\mathbf C}P} 
\newcommand{\bstar}{{\mathbf C}^*} 
\newcommand{\cinf}{C_\infty} 
\newcommand{\zz}{{\mathbf Z}}

\newcommand{\minus}{{-1}}

\newcommand{\eltwo}{{\ell^2}}

\newcommand{\bZ}{{\mathbf Z}}

\newcommand{\bC}{{\mathbf C}}

\newcommand{\wt}{\widetilde}

\newcommand{\ol}{\overline}

\newcommand{\Min}{\mathrm{Min}}
\newcommand{\MU}{\Min (U)} 
\newcommand{\rk}{\mathrm{rk}}
\newcommand{\GR}{\operatorname{Gr}} 
\newcommand{\grh}{\GR{H}}

\newcommand{\sg}{\mathrm {sing}} 
 
\def\clap#1{\hbox to 0pt{\hss#1\hss}}

\newcommand{\comment}[1]{} 
 
\newcommand{\gb}{\beta}

\newcommand{\gr}{\rho} 
\newcommand{\gs}{\sigma}

\newcommand{\gS}{\Sigma} 

\newcommand{\OL}{{\ol{L}}}

\newcommand{\Hom}{\operatorname{Hom}}

	{ 
\end{enumerate}
} 
\newenvironment{enumerate1*}{ 
\begin{enumerate}[\upshape (*1)]}
	{ 
\end{enumerate}
} 
	{ 
\end{enumerate}
}

\newenvironment{enumeratea'}{ 
\begin{enumerate}[\upshape (a)$'$]}{ 
\end{enumerate}
} 

\begin{document} 

\title{Cohomology of hyperplane complements with group ring coefficients} 
\author{Michael W. Davis\thanks{The first and the second authors were partially
    supported by NSF grant DMS 0706259.}\and{Tadeusz Januszkiewicz$^*$}
		\and{Ian J. Leary}\thanks{The
    third author was partially supported by NSF grant DMS 0804226.}
  \and{Boris Okun}}

\date{\today} \maketitle 
\begin{abstract}
	We compute the cohomology with group ring coefficients of the complement of a finite collection of affine hyperplanes in $\bC^n$.
	It is nonzero in exactly one degree, namely, the degree equal to the rank of the arrangement.


	\end{abstract}

\section*{ }
A \emph{hyperplane arrangement} $\ca$ is a finite collection of affine hyperplanes in $\bC^n$.
A \emph{subspace} of $\ca$ is a nonempty intersection of hyperplanes in $\ca$.
Denote by $L(\ca)$ the poset of subspaces, ordered by inclusion.
Put $\ol{L}(\ca):=L(\ca)\cup\{\bC^n\}$.
An arrangement is \emph{central} if $L(\ca)$ has a unique minimum element.
In general, the minimal elements of $L(\ca)$ are a family of parallel subspaces.
The \emph{rank} of $\ca$ is the codimension in $\bC^n$ of a minimal element.
$\ca$ is \emph{essential} if $\rk(\ca)=n$.
Given $G\in \OL(\ca)$, put 
\[
	\ca_G:=\{H\in \ca\mid H\supseteq G\}.
\]
It is a central arrangement of rank $\gr(G)=n-d(G)$, where $d(G)=\dim_\bC G$.

The \emph{singular set} $\gS(\ca)$ of the arrangement is the union of hyperplanes in $\ca$ (so, $\gS(\ca)$ is a subset of $\bC^n$).
The complement of $\gS(\ca)$ in $\bC^n$ is denoted $M(\ca)$.
Similarly, the complement of $\gS(\ca_G)$ in $\bC^n$ is $M(\ca_G)$.

We now state our main result.  

\begin{theorem}\label{t:main}
	Suppose $\ca$ is an arrangement of rank $l$.
	Let $\pi=\pi_1(M(\ca))$.
	Then $H^*(M(\ca);\zz \pi)$ is concentrated in degree $l$ and is free abelian.
\end{theorem}
\begin{corollary}\label{c:FL}
	The right $\zz\pi$-module $H^l(M(\ca);\zz\pi)$ is type $F\!L$.
\end{corollary}
\begin{proof}
	It is known that $M(\ca)$ is homotopy equivalent to a finite
	complex $X$ of dimension $l$. 
	Using the cellular cochains of the universal cover of $X$ we
	get a free resolution of length $l$:  
	\[
		0\to C^0(X;\zz\pi)\to \cdots \to C^{l}(X;\zz\pi)\to
		H^l(X;\zz\pi)\to 0. 
	\]
\end{proof}

A group $\pi$ is a \emph{duality group} if it is type $F\!P$ and
$H^*(\pi;\zz \pi)$ is concentrated in a single degree and is
torsion-free. 
\begin{corollary}[cf.~\cite{squirer}]\label{c:main}
	Suppose $\ca$ is a $K(\pi,1)$ arrangement (i.e., $M(\ca)$ is a
	$K(\pi,1)$ with $\pi=\pi_1(M(\ca))$). 
	Then $\pi$ is a duality group.
\end{corollary}

The next lemma is well-known.
\begin{lemma}[{cf.~\cite[Prop.~2.1]{djl}}]\label{l:sigma}
	Suppose $\ca$ is a hyperplane arrangement of rank $l$.
	Then $\gS(\ca)$ is homotopy equivalent to a wedge of $(l-1)$-spheres.
\end{lemma}

For each $G\in \OL(\ca)$, $\ca\cap G$ denotes the hyperplane arrangement in $G$ consisting of all elements of $L(\ca)$ which are subspaces of codimension-one in $G$.
Then $\ca\cap G$ is an arrangement of rank $l(G)=d(G)-n_0$, where $n_0$ is the rank of a minimal element of $L(\ca)$.
We note that 
\begin{equation}\label{e:rank}
	l(G)+\gr(G)=n-n_0=l.
\end{equation}

Let $\gb(\ca\cap G)$ denote the reduced Betti number of $G\cap\gS$ in degree $l(G)-1$, i.e., 
\begin{equation}\label{e:betti}
	\gb(\ca\cap G):=\rk (H^{l(G)}(G,\gS(\ca\cap G))).
\end{equation}

Suppose $\ca$ is an essential, central arrangement in $\bC^n$.
Projectivizing we get a projective hyperplane arrangement in $P\ca$ in $\bc^{n-1}$.
Choose a hyperplane in $P\ca$ to regard as the hyperplane at infinity.
Removing it, we obtain a hyperplane arrangement $\ca'$ in $\bC^{n-1}$, called an \emph{associated affine arrangement}.
We note that $M(\ca)$ is a $\bstar$-bundle over $M(\ca')$; moreover, this bundle is trivial (since either $n=1$ or $\ca'$ is nonempty).
Thus, $M(\ca)\cong M(\ca')\times \bstar$.
Let $\cinf$ denote the fundamental group of $\bstar$ (i.e., $\cinf$ is
the infinite cyclic group).  From the above discussion we get the
following. 
\begin{lemma}\label{l:product}
	Suppose $\ca$ is an essential, central arrangement in $\bC^n$ and $\ca'$ is an associated affine arrangement.
	Put $\pi=\pi_1(M(\ca))$, $\pi'=\pi_1(M(\ca'))$.
	Then $\pi=\pi'\times \cinf$, and 
	\[
		H^*(M(\ca); \zz\pi)=H^{*-1}(M(\ca');\zz\pi')\otimes \zz, 
	\]
	where $\cinf$ acts trivially on $\bZ$.
\end{lemma}
\begin{proof}
	\[
		H^i(\bstar;\zz\cinf)=H^i(S^1;\zz\cinf)= 
		\begin{cases}
			\zz, &\text{if $i=1$;}\\
			0,&\text{if $i\neq 1$.} 
		\end{cases}
	\]
	So, the equation in the lemma follows from the K\"unneth Formula.
\end{proof}

Suppose $\ca$ is a hyperplane arrangement in $\bC^n$.
An open convex subset $U$ in $\bC^n$ is \emph{small} (with respect to $\ca$) if $\{G\in \ol{L}(\ca)\mid G\cap U\neq\emptyset\}$ has a unique minimum element $\MU$.  The intersection of two small convex open sets is also small; hence, the same is true for any finite intersection of such sets.

Now let $\cu=\{U_i\}_{i\in I}$ be an open cover of $\bC^n$ by small convex sets.
We may suppose that $\cu$ is finite and that it is closed under taking intersections.
For each $G\in \ol{L}(\ca)$, put 
\begin{align*}
	\cu_G:\!\!&=\{U\in \cu\mid \MU\subseteq G\},\\
	\cu^\sg_G:\!\!&=\{U\in \cu\mid \MU\subsetneq G\}=\{U\in
	\cu_G\mid U\cap \gS(\ca \cap G)\neq\emptyset\}.\\
\end{align*}
The open cover $\cu$ restricts to an open cover $\widehat{\cu}=\{U-\gS(\ca)\}_{U\in \cu}$ of $M(\ca)$.
Any element $\widehat{U}=U-\gS(\ca)$ of the cover is homotopy equivalent to the complement of a central arrangement $M(\ca_G)$, where $G=\MU$.

Suppose $N(\cu)$ is the nerve of $\cu$ and $N(\cu_G)$ is the subcomplex defined by $\cu_G$.
Since $N(\cu_G)$ and $N(\cu^\sg_G)$ are nerves of covers of $G$ and $\gS(\ca\cap G)$, respectively, by contractible open subsets, we have that for each $G\in \OL(\ca)$, 
\begin{equation}\label{e:nerve}
	H^*(N(\cu_G),N(\cu^\sg_G))=H^*(G,\gS(\ca\cap G)).
\end{equation}
For each $k$-simplex $\gs=\{i_0,\dots, i_k\}$ in $N(\cu)$, let 
\[
	U_\gs:=U_{i_0}\cap \cdots \cap U_{i_k} 
\]
denote the corresponding intersection.

Let $r:\wt{M}(\ca)\to M(\ca)$ be the universal cover.
The induced cover $\{r^\minus(\widehat{U})\}$ of $\wt{M}(\ca)$ has the same nerve $N(\widehat{\cu})$ ($=N(\cu)$).
We have the Mayer--Vietoris double complex, 
\[
	C_{i,j}:=\bigoplus_{\gs\in N^{(i)}} C_j(r^\minus (\widehat{U}_\gs)), 
\]
where $N^{(i)}$ denotes the set of $i$-simplices in $N(\cu)$ (cf.~\cite[Ch.~VII]{brown}.) We get a corresponding double cochain complex, 
\begin{equation}\label{e:e0}
	E^{i,j}_0:= \Hom_\pi(C_{i,j},\zz \pi), 
\end{equation}
where $\pi= \pi_1(M(\ca))$.
The filtration on the double complex gives a spectral sequence converging to the associated graded module for cohomology: 
\[
	\grh^m(M(\ca);\zz\pi)=E_\infty := \bigoplus_{i+j=m} E^{i,j}_\infty . 
\]
\begin{proof}[Proof of Theorem~\ref{t:main}] The proof is by induction on the rank $l$ of $\ca$.
	The result is trivial for $l=0$ (for then the arrangement is empty).
	Lemma~\ref{l:product} shows that if we know the result for ranks less than $l$, then we also know it for any central arrangement of rank $l$.
	So, given a rank $l$ arrangement $\ca$, the inductive hypothesis implies that the theorem holds for each small open set in our cover $\cu$.
	In other words, we can assume that for each $U\in \cu$, for $G=\MU$ and $\pi_G=\pi_1(M(\ca_G))$, $H^*(U-\gS;\zz\pi_G)=H^*(M(\ca_G);\zz\pi_G)$ is free abelian and is concentrated in degree $\gr(G)=l-l(G)$.
	
	By first using the horizontal differential in \eqref{e:e0}, we get a spectral sequence with $E_1$-terms 
	\begin{equation}\label{e:e1}
		E^{i,j}_1= C^i(N(\cu);\ch^j), 
	\end{equation}
	where $\ch^j$ is the coefficient system on $N(\cu)$ defined by 
	\[
		\gs\mapsto H^j(M(\ca_G);\zz\pi), 
	\]
	for $G=\Min(U_\gs)$.
	These coefficients are $0$ for $j\neq \gr(G)$, i.e., for $l(G)\neq l-j$ (by \eqref{e:rank}).
	Moreover, for any coface $\gs'$ of $\gs$, if $G':=\Min(U_{\gs'})\subsetneq G$, then the coefficient homomorphism $H^j(M(\ca_G);\zz \pi)\to H^j(M(\ca_{G'});\zz\pi)$ is the zero map.
	It follows that the $E_1$ page of the spectral sequence decomposes as a direct sum (cf.~\cite[Lemma~2.2]{do}).
	For a fixed $j$, the $E^{i,j}_1$ term decomposes as 
	\begin{equation*}
		E^{i,j}_1=\bigoplus_{G\in \OL_{n-j}(\ca)} C^i(N(\cu_G),N(\cu^\sg_G);H^j(M(\ca_G);\zz \pi)), 
	\end{equation*}
	where we have constant coefficients in each summand.
	Hence, at $E_2$ we have 
	\begin{align*}
		E^{i,j}_2&=\bigoplus_{G\in \OL_{n-j}(\ca)} H^i(N(\cu_G),N(\cu^\sg_G);H^j(M(\ca_G);\zz \pi))\\
		&=\bigoplus_{G\in \OL_{n-j}(\ca)} H^i(G,\gS(\ca\cap G);H^j(M(\ca_G);\zz \pi)),
	\end{align*}
	where the second equation is by \eqref{e:nerve}.
	By Lemma~\ref{l:sigma}, $H^i(G,\gS(\ca\cap G))$ is nonzero only for $i=l(G)=l-j$.
	So, the $E_2$ terms are nonzero only in total degree $l$.
	It follows that the spectral sequence collapses at $E_2$.
	Thus, for $k\neq l$, $H^k(M(\ca);\zz\pi)=0$, while 
	\begin{equation}\label{e:main}
		\grh^l(M(\ca);\zz\pi)=\bigoplus_{G\in \OL(\ca)} H^{l(G)}(G,\gS(G\cap \ca);H^{\gr(G)}(M(\ca_G);\zz\pi)),
	\end{equation}
	and, therefore, is free abelian. It follows that 
	$H^l(M(\ca);\zz\pi)$, the ungraded object, is also free abelian. 
\end{proof}
\begin{remark}\rm 
	Here are some more comments about \eqref{e:main}.
	Since $\zz \pi$ is a $\pi$-bimodule, $H^l(M(\ca);\zz\pi)$ is a right $\zz \pi$-module and $\grh^l(M(\ca);\zz\pi)$ is the associated graded $\zz\pi$-module.
	Similarly, each summand on the right hand side of \eqref{e:main} is a $\zz\pi$-module and the formula is an isomorphism of $\zz\pi$-modules.
	
	The coefficients in the summand corresponding to $G$ come from the induced representation, 
	\[
		H^{\gr(G)}(M(\ca_G);\zz\pi)=H^{\gr(G)}(M(\ca_G);\zz\pi_G)\otimes _{\pi_G} \zz \pi, 
	\]
	where $\pi_G:=\pi_1(M(\ca_G))$.
	So, the summand corresponding to $G$ is a sum of $\gb(\ca \cap G)$ copies of the induced representation \( H^{\gr(G)}(M(\ca_G);\zz\pi_G)\otimes _{\pi_G} \zz \pi, \) where $\gb(\ca\cap G)$ was defined in \eqref{e:betti}.
	If $\ca_G\neq \emptyset$, $H^{\gr(G)}(M(\ca_G);\zz\pi_G)$ is not a free $\zz\pi_G$-module.
	The reason is that if $M(\ca'_G)$ is an associated affine arrangement to $\ca_G$ and $\pi'_G=\pi_1(\ca'_G)$, then, by Lemma~\ref{l:product}, 
	\[
		H^{\gr(G)}(M(\ca_G);\zz\pi_G)=H^{\gr(G)-1}(M(\ca'_G);\zz\pi'_G)\otimes \zz, 
	\]
	which is not free (unless $\pi_G =1$).
	Hence, only one summand on the right hand side of \eqref{e:main} is a free $\zz\pi$-module, the one corresponding to $G=\bC^n$.
	It is 
	\[
		H^l(\bC^n,\gS(\ca)) \otimes \zz\pi, 
	\]
	which is a free of rank $\gb(\ca)$.
	In \cite[Theorem 6.2]{djl} we showed that the reduced $\ell^2$-cohomology of $M(\ca)$ is $H^l(\bC^n,\gS(\ca)) \otimes \ell^2\pi$.
	The free summand described above injects into $\ell^2$-cohomology, while the other summands map to $0$.
\end{remark}


\bigskip
\leftline{\bf Addresses:} 

\smallskip
\noindent 
Michael W. Davis, Tadeusz Januszkiewicz and Ian J. Leary

\noindent
Department of Mathematics, The Ohio State University, 231 W. 18th Ave., Columbus OH 43210, USA.  

\smallskip
\noindent
Tadeusz Januszkiewicz 

\noindent
The Mathematical Institute of the Polish Academy of Sciences.  

\noindent 
On leave from Instytut Matematyczny, Uniwersytet Wroc\l awski.  

\smallskip
\noindent 
Boris Okun 

\noindent 
Department of Mathematics, University of Wisconsin-Milwaukee, Milwaukee, WI 53201-0413, USA

\smallskip
\leftline{\bf Email:} 

\noindent 
{\tt mdavis@math.ohio-state.edu \qquad 
tjan@math.ohio-state.edu 

\noindent 
leary@math.ohio-state.edu\phantom{m} \qquad
okun@uwm.edu}

\end{document}